\documentclass[twoside,12pt]{article}

\usepackage[T1]{fontenc}
\usepackage{lmodern}
\usepackage{amsmath,amssymb,amsthm,mathtools,mathrsfs}
\usepackage{graphicx}
\usepackage{xcolor}
\usepackage{microtype}
\usepackage{enumitem}
\usepackage[top=2cm,bottom=2.5cm,left=2.5cm,right=2.5cm]{geometry}
\usepackage{booktabs}

\newcommand{\floor}[1]{\left\lfloor #1\right\rfloor}

\numberwithin{equation}{section}

\theoremstyle{plain}
\newtheorem{theorem}{Theorem}[section]
\newtheorem{proposition}[theorem]{Proposition}
\newtheorem{lemma}[theorem]{Lemma}
\newtheorem{corollary}[theorem]{Corollary}
\theoremstyle{definition}
\newtheorem{definition}[theorem]{Definition}
\theoremstyle{remark}

\usepackage[colorlinks=true,linkcolor=blue!55!black,citecolor=blue!55!black,urlcolor=blue!55!black]{hyperref}
\usepackage[nameinlink,capitalize,noabbrev]{cleveref}
\crefname{theorem}{Theorem}{Theorems}
\crefname{proposition}{Proposition}{Propositions}
\crefname{lemma}{Lemma}{Lemmas}
\crefname{corollary}{Corollary}{Corollaries}
\crefname{definition}{Definition}{Definitions}
\crefname{remark}{Remark}{Remarks}
\Crefname{theorem}{Theorem}{Theorems}
\Crefname{proposition}{Proposition}{Propositions}
\Crefname{lemma}{Lemma}{Lemmas}
\Crefname{corollary}{Corollary}{Corollaries}
\Crefname{definition}{Definition}{Definitions}
\Crefname{remark}{Remark}{Remarks}

\newcommand{\N}{\mathbb N}
\newcommand{\cA}{\mathcal A}
\newcommand{\cB}{\mathcal B}
\newcommand{\cD}{\mathcal D}
\newcommand{\cE}{\mathcal E}
\newcommand{\cN}{\mathcal N}
\newcommand{\Pw}{P_{\mathrm w}}
\newcommand{\op}{\mathbin{\oplus}}
\newcommand{\red}{\operatorname{red}}
\newcommand{\rev}{\operatorname{rev}}
\newcommand{\dir}{\operatorname{dir}}

\hypersetup{
  pdftitle={Short Permutation Pattern Waves: The Complete Length-Five Classification},
  pdfauthor={Mengmeng He, Ralf Klasing, Yaping Mao, and Xiaoyan Zhang}
}

\title{\textbf{Short permutation pattern waves: The complete length-five classification}%
\thanks{Supported by the National Natural Science Foundation of China
(Nos.~12271259, 12471329, 12061059 and 12201407).}}

\author{Lihong Qiu\thanks{School of Mathematics, Faculty of Basic Sciences, University of Shanghai for Science and Technology, Shanghai 200093, China. \texttt{qiulihong@usst.edu.cn}}, \ \ Minglei Mao\thanks{School of
Mathematics and Statistics, Qinghai Normal University, Xining, Qinghai
810008, China. \texttt{mmldyxa@163.com}}, \ \ Yaping Mao%
\thanks{Academy of Plateau Science and Sustainability, and School of
Mathematics and Statistics, Qinghai Normal University, Xining, Qinghai
810008, China. \texttt{yapingmao@outlook.com; myp@qhnu.edu.cn}}}
\date{}

\begin{document}

\maketitle

\begin{abstract}
For a permutation $\pi\in S_k$, a $\pi$-wave is an increasing sequence
of integers whose consecutive gaps have relative order $\pi$. Let $g(\pi,n)$
be the maximum cardinality of a $\pi$-wave-free subset of $[n]$, and let $P(\pi,r)$
be the least integer $N$ such that every $r$-coloring of $[N]$ contains a monochromatic
$\pi$-wave. We determine the asymptotic orders of both functions for every
permutation of length at most five. Existing upper bounds and weak-wave constructions provide sharp bounds for all  but eight length-five patterns in $S_{5}$. 
For these eight exceptional alternating patterns, we establish \(g(\pi,n)=\Theta((\log n)^2)\) and \(P(\pi,r)=\Theta(r(\log r)^2)\).\\[0.2mm]
\end{abstract}

\section{Introduction}

Throughout this paper,
 we write $\N=\{1,2,3,\dots\}$ and $[n]=\{1,2,\dots,n\}$.
Let $S_k$ denote the symmetric group on $[k]$. We write permutations in one-line notation. For a word $w=w_1\cdots w_m$ of distinct real numbers, its
\emph{normalization} or \emph{reduction} is the unique permutation
$\red(w)\in S_m$ such that $\red(w)(i)<\red(w)(j)$ if and only if
$w_i<w_j$. Thus $\red(w)$ records the relative order of the entries of $w$.
We apply this operation to the consecutive gaps of an increasing sequence.

\begin{definition}\label{def:pi-wave}
Let $\pi\in S_k$.  An increasing sequence
$x_1<\cdots<x_{k+1}$ of integers is a \emph{$\pi$-wave} if its consecutive gaps are
pairwise distinct and $\red(x_2-x_1,x_3-x_2,\dots,x_{k+1}-x_k)=\pi.$
Equivalently, for all $i,j\in[k]$,
\[
 x_{i+1}-x_i>x_{j+1}-x_j
 \quad\Longleftrightarrow\quad
 \pi(i)>\pi(j).
\]
A set of integers is \emph{$\pi$-wave-free} if it contains no $\pi$-wave.
\end{definition}

For $\pi\in S_k$, let $ g(\pi,n)=\max\{|A|:A\subseteq[n]\text{ is $\pi$-wave-free}\}$. An \emph{$r$-coloring} of a set $X$ is a map $c:X\to[r]$, and a subset $Y\subseteq
X$ is \emph{monochromatic} if $c$ is constant on $Y$. Let $P(\pi,r)$ be the
least positive integer $N$ such that every $r$-coloring of $[N]$ contains a
monochromatic $\pi$-wave. The existence of $P(\pi,r)$ follows from van der Waerden's theorem \cite{VanDerWaerden1927}. 
Indeed, let
$ L=1+\sum_{i=1}^k\pi(i)=1+\frac{k(k+1)}2$.  Every sufficiently long $r$-colored interval
contains a monochromatic arithmetic progression
$a,a+d,\dots,a+(L-1)d$, where $d>0$. The points \(y_1=a\) and \(y_{i+1}=a+d\sum_{j=1}^i\pi(j)\) for \(1\le i\le k\) lie in this progression and 
have consecutive gaps \(d\pi(1),\dots,d\pi(k)\); hence they form a monochromatic \(\pi\)-wave.

We use the following standing conventions.  All asymptotic assertions are
taken as $n\to\infty$ or $r\to\infty$, with the permutation fixed.  Unless
otherwise indicated, the constants in $O$, $\Omega$, and
$\Theta$ may depend on the fixed permutation. A subscript in
$O_\pi(\cdot)$ makes this dependence explicit.
All logarithms are natural unless a base is indicated. Unless explicitly
qualified as weak, a wave always means the strict notion in
Definition \ref{def:pi-wave}.  The classification depends on a few structural features of permutations. An index $i\in\{2,\dots,k-1\}$ is a \emph{peak} of $\pi$ if
$\pi(i)>\pi(i-1)$ and $\pi(i)>\pi(i+1)$.  For a set $X$ of values occurring in $\pi$,
let $\pi\setminus X$ be the normalization of the word obtained by
deleting those values.  We use 
$\pi^{\rev}=\pi(k)\pi(k-1)\cdots\pi(1)$ for the reversal of $\pi$.

The study of waves began with monotone sequences of consecutive
gaps. Brown, Erd\H{o}s, and Freedman~\cite{BrownErdosFreedman} introduced
descending waves, allowing equal consecutive gaps. Alon and
Spencer~\cite{AlonSpencer} determined the cubic order of the two-color
threshold for these weak waves as the wave length tends to infinity. Reflection of the ambient
interval gives the same threshold for ascending waves.
LeSaulnier and Robertson~\cite{LeSaulnierRobertson} extended the analysis
to more colors. Related variants include zero-sum ascending waves, studied by
A.~Bialostocki, G.~Bialostocki, Caro, and Yuster~\cite{BialostockiEtAl},
and intermingled ascending wave sets, studied by Robertson, Cremin, Daniel,
and Xiang~\cite{CreminEtAl}. In a recent preprint, Mao~\cite{MaoDense2026}
presents a resolution of the Alon--Spencer conjecture on the length of weak
ascending waves in dense subsets of $[n]$.
A different line of work fixes the wave length and lets the number
of colors grow. Bollob\'as, Erd\H{o}s, and Jin~\cite{BollobasErdosJin}
proved that, for a monotone gap pattern of fixed length $k$, the
monochromatic threshold has order $r(\log r)^{k-1}$.
Landman and Robertson~\cite{LandmanRobertson} established the quartic order for the two-color threshold
of strictly monotone waves as their length tends to infinity,  and extended the setting to arbitrary permutations. Cong~\cite{Cong} approached this question through the extremal
function $g(\pi,n)$. He proved the general bound
$g(\pi,n)=O_\pi((\log n)^{k-1})$ for $\pi\in S_k$ and showed that this
order is attained precisely when $\pi$ has no peak. 

Some other Ramsey problems restrict the consecutive gaps to prescribed sets.
Landman, Robertson, and Robertson~\cite{LandmanRobertsonRobertson2023}
studied accessibility and Ramsey functions for such restricted-gap
sequences. Clifton~\cite{Clifton2024} obtained new bounds for
diffsequences, and Wesley~\cite{Wesley2025} investigated corresponding
questions for Fibonacci gaps. Quester~\cite{Quester2025} studied structural
properties of accessible sets.

The density--coloring lemma of Cong~\cite{Cong}
links the two parameters: an upper bound $g(\pi,n)=O((\log n)^t)$ implies
$P(\pi,r)=O(r(\log r)^t)$, while a lower bound
$P(\pi,r)=\Omega(r(\log r)^t)$ implies
$g(\pi,n)=\Omega((\log n)^t)$. These results allow us to combine
extremal upper bounds with coloring lower bounds.
Cong~\cite{Cong} asked whether every fixed permutation
admits an exponent $e_\pi$ such that
$g(\pi,n)=\Theta((\log n)^{e_\pi})$.

 We determine the asymptotic orders of both functions for every permutation of length at most five. 
Applying the existing upper bounds and weak-wave lower-bound constructions yields tight bounds for all patterns of length at most four and for $112$ of the $120$ length-five patterns. 
We give the resulting finite classification explicitly. For the remaining eight up--down patterns, these bounds leave a gap between $\Omega((\log n)^2)$ and $O((\log n)^3)$. 
Our main contribution is a quadratic logarithmic upper bound for these eight patterns, completing the classification.
 
We state the classification in increasing order of pattern length,
beginning with the patterns of length at most three.

\begin{theorem}\label{thm:small}
The complete classification for $k\le3$ is as follows.
\begin{enumerate}[label=\textup{(\roman*)},leftmargin=2.1em]
\item For the unique pattern $1\in S_1$,
 $g(1,n)=1$,
 $P(1,r)=r+1$.
\item For every $\pi\in S_2$,
 $g(\pi,n)=\Theta(\log n)$, 
 $P(\pi,r)=\Theta(r\log r)$.
\item If $\pi\in S_3$ has no peak, equivalently
$\pi\in\{123,213,312,321\}$, then
$g(\pi,n)=\Theta((\log n)^2)$ and $P(\pi,r)=\Theta(r(\log r)^2).$
For the two remaining patterns $\pi\in\{132,231\}$,
$
 g(\pi,n)=\Theta(\log n),
 P(\pi,r)=\Theta(r\log r).
$
\end{enumerate}
\end{theorem}
For length four, the same distinction gives the complete classification.
\begin{theorem}\label{thm:S4}
Let $\pi\in S_4$.  Then
\[
 g(\pi,n)=
 \begin{cases}
  \Theta((\log n)^3),&\pi\text{ has no peak},\\
  \Theta((\log n)^2),&\pi\text{ has a peak},
 \end{cases}
\]
and the corresponding coloring orders are
\[
 P(\pi,r)=
 \begin{cases}
  \Theta(r(\log r)^3),&\pi\text{ has no peak},\\
  \Theta(r(\log r)^2),&\pi\text{ has a peak}.
 \end{cases}
\]
\end{theorem}

For $\pi\in S_5$, the presence of a peak no longer determines the exponents. Define the \emph{two-value deletion class}
\[
 \cE_2=\{\pi\in S_5:1,2\text{ are nonadjacent and }
                    \pi\setminus\{1,2\}\in\{132,231\}\}
\]
and the \emph{exceptional alternating class}
\[
 \cA_5=\{14352,15243,15342,24351,25143,25341,34152,34251\}.
\]
Every member of $\cA_5$ is up--down alternating, meaning
$\pi(1)<\pi(2)>\pi(3)<\pi(4)>\pi(5)$. We call these patterns ``exceptional'' 
because they are precisely the eight patterns left unresolved by the deletion bounds
 and the lower-bound constructions used below. The class $\cA_5$ is a proper subset of the up-down
 alternating permutations of length five. The complete length-five classification is as follows.

\begin{theorem}\label{thm:S5}
Let $\pi\in S_5$.
\begin{enumerate}[label=\textup{(\roman*)},leftmargin=2.1em]
\item If $\pi$ has no peak, then
$
 g(\pi,n)=\Theta((\log n)^4),
 P(\pi,r)=\Theta(r(\log r)^4).
$
There are $16$ such patterns.
\item If $\pi\in\cE_2$, then
$g(\pi,n)=\Theta((\log n)^2),
 P(\pi,r)=\Theta(r(\log r)^2).$
There are $24$ such patterns.
\item If $\pi$ has a peak and $\pi\notin\cE_2\cup\cA_5$, then
$ g(\pi,n)=\Theta((\log n)^3),P(\pi,r)=\Theta(r(\log r)^3).$
There are $72$ such patterns.
\item If $\pi\in\cA_5$, then
$ g(\pi,n)=\Theta((\log n)^2),
 P(\pi,r)=\Theta(r(\log r)^2).$
\end{enumerate}
\end{theorem}

The rest of the paper is arranged as follows. We give some preliminaries and lower-bound constructions in Section~\ref{sec:lower}.
In Section~\ref{sec:classification}, we establish the shorter classifications and
reduces the length-five case to an upper bound for $\cA_5$.
Section~\ref{sec:separator} develops the
separator-tree construction and bounds for labelled binary trees in terms
of child-subtree label overlaps. In Section~\ref{sec:completion}, we combine these
bounds with local completion and scale--depth counting to prove the remaining
upper bound and complete Theorem~\ref{thm:S5}. 
Section~\ref{sec:conclusion} discusses a direction for further research.

\section{Preliminaries and lower-bound constructions}\label{sec:lower}
For lower bounds we use the weak-difference waves introduced by Cong \cite[Definition~4.2]{Cong}.  A \emph{weak-difference
$\pi$-wave}, abbreviated to a \emph{weak $\pi$-wave}, is an increasing
sequence $x_1<\cdots<x_{k+1}$ such that $\pi(i)>\pi(j)\Longrightarrow\
 x_{i+1}-x_i\ge x_{j+1}-x_j.$ Let $\Pw(\pi,r)$ be the least $N$ forcing a monochromatic weak $\pi$-wave in
every $r$-coloring of $[N]$.  Since every strict $\pi$-wave is weak,
 $\Pw(\pi,r)\le P(\pi,r)$. For $1\le a\le b\le k$, the consecutive pattern $\pi[a,b]$ is 
 the normalization $\red(\pi(a),\pi(a+1),\dots,\pi(b))$.

\begin{lemma}\label{lem:consecutive}
If $\sigma$ is a consecutive pattern of $\pi$, then
$g(\pi,n)\ge g(\sigma,n)$ and $P(\pi,r)\ge P(\sigma,r).$
The same assertion holds with $P$ replaced by $\Pw$.
\end{lemma}

\begin{proof}
Suppose that $\sigma=\pi[a,b]$.  If
$x_1<\cdots<x_{k+1}$ is a $\pi$-wave, then
$x_a<x_{a+1}<\cdots<x_{b+1}$ is a $\sigma$-wave.  Hence every
$\sigma$-wave-free set is $\pi$-wave-free, and every coloring with no
monochromatic $\sigma$-wave has no monochromatic $\pi$-wave.  The same
argument applies to weak waves, since passing to 
consecutive gaps preserves the required inequalities.
\end{proof}

\begin{lemma}\label{lem:reversal-invariance}
Let $\pi^{\rev}$ denote the reversal of the word $\pi$.  Then
\[
 g(\pi,n)=g(\pi^{\rev},n),\qquad
 P(\pi,r)=P(\pi^{\rev},r),\qquad
 \Pw(\pi,r)=\Pw(\pi^{\rev},r).
\]
\end{lemma}

\begin{proof}
Reflect $[n]$ by $x\mapsto n+1-x$ and reorder each selected
sequence increasingly. This reverses the order of the consecutive gaps and preserves all
comparisons between them. Applying the same reflection to each color class
proves the two coloring identities, including the weak-wave version.
\end{proof}

The following density-coloring related lemma is due to~\cite{Cong}.

\begin{lemma} \cite{Cong} \label{lem:transfer}
For fixed $\pi$ and fixed $t\ge0$:
\begin{enumerate}[label=\textup{(\roman*)},leftmargin=2.1em]
\item If $g(\pi,n)=O((\log n)^t)$, then
$P(\pi,r)=O(r(\log r)^t)$.
\item If $P(\pi,r)=\Omega(r(\log r)^t)$, then
$g(\pi,n)=\Omega((\log n)^t)$.
\end{enumerate}
Consequently, a lower bound of order $r(\log r)^t$ for $\Pw(\pi,r)$ gives the
same exponent as a lower bound for $g(\pi,n)$.
\end{lemma}

The following endpoint-extension bound for weak waves will be used in the subsequent
lower-bound constructions.

\begin{lemma}[\cite{Cong}]\label{lem:endpoint} Let $k\geq 2$ and $t\geq 0$.
Suppose that $\pi\in S_k$ begins with $k$, and let $\sigma$ be obtained by
deleting this first entry.  If
$\Pw(\sigma,r)=\Omega(r(\log r)^t),$
then
$\Pw(\pi,r)=\Omega(r(\log r)^{t+1}).$
The same conclusion holds when $\pi$ ends with $k$.
\end{lemma}

\begin{proof}
When $\pi$ begins with its maximum entry, the assertion follows
from Cong~\cite[Proposition~1.10 and Remark~4.3]{Cong}.
When $\pi$ ends with its maximum entry, apply the first case
to $\pi^{\rev}$ and use Lemma~\ref{lem:reversal-invariance}.
\end{proof}

In a permutation with no peak, the maximum occurs at an endpoint. 
Applying Lemma \ref{lem:endpoint} repeatedly yields the following lower bound.

\begin{corollary}\label{cor:nopeakweak}
If $\pi\in S_k$ has no peak, then
$\Pw(\pi,r)=\Omega(r(\log r)^{k-1}).$
Consequently,
$g(\pi,n)=\Omega((\log n)^{k-1})$.
\end{corollary}

\begin{proof}
We proceed by induction on $k$.  For $k=1$, every color class containing two points gives
a weak one-gap wave, so $\Pw(1,r)=r+1=\Omega(r)$.  If $k\ge2$, the maximum
value $k$ must occur at one of the two ends of a no-peak permutation;
otherwise its position would be a peak.  Deleting that endpoint leaves a
no-peak permutation of length $k-1$.  The induction hypothesis and
Lemma \ref{lem:endpoint} yield one additional factor of $\log r$.  The density lower bound then
follows from Lemma \ref{lem:transfer}.
\end{proof}

We next give the weak-wave product construction of
Cong~\cite[Lemma~4.5]{Cong} for skew sums, and obtain the direct-sum version
by reversal. The proof below uses integer coordinates throughout.  For
$\alpha\in S_a$ and $\beta\in S_b$, their direct and skew sums are
$\alpha\op\beta
 =\alpha(1)\cdots\alpha(a)(\beta(1)+a)\cdots(\beta(b)+a)$ and
$\alpha\ominus\beta
 =(\alpha(1)+b)\cdots(\alpha(a)+b)\beta(1)\cdots\beta(b)$. Thus the gaps prescribed by $\alpha$ are smaller than those prescribed by
$\beta$ in $\alpha\op\beta$, and larger in $\alpha\ominus\beta$.

\begin{proposition}\label{prop:product}
Let $\alpha\in S_a$ and $\beta\in S_b$.  Suppose $p\geq 0$ and $q\geq0$. If
$\Pw(\alpha,r)=\Omega(r(\log r)^p)$ and
 $\Pw(\beta,r)=\Omega(r(\log r)^q),$
then 
\[
\Pw(\alpha\ominus\beta,r)
 =\Omega(r(\log r)^{p+q})
 \quad\text{and}\quad  
 \Pw(\alpha\oplus\beta,r)
 =\Omega(r(\log r)^{p+q}).
\]
\end{proposition}

\begin{proof}
First we prove the skew-sum assertion.  Then the direct-sum assertion follows by
applying the skew-sum result to $\beta^{\rev}$ and $\alpha^{\rev}$, because $(\alpha\oplus\beta)^{\rev}
 =\beta^{\rev}\ominus\alpha^{\rev}$, and the two hypotheses are invariant under reversal by
Lemma \ref{lem:reversal-invariance}.

Fix large $r$ and let $m=\lfloor\sqrt{r/5}\rfloor$.  Let
$L_\alpha=\Pw(\alpha,m)-1,
 L_\beta=\Pw(\beta,m)-1$.
Choose $m$-colorings $c_\alpha:[L_\alpha]\to[m]$ and
$c_\beta:[L_\beta]\to[m]$ containing no monochromatic weak $\alpha$-wave
and no monochromatic weak $\beta$-wave, respectively.  Every integer
$x\in[5L_\alpha L_\beta]$ has a unique representation
$x=5L_\beta(u-1)+L_\beta(h-1)+v,
 u\in[L_\alpha],\ h\in[5],\ v\in[L_\beta]$.
Call $u,h,v$ the coarse, slab, and fine coordinates of $x$, respectively.
Assign to $x$ the color
 $(c_\alpha(u),c_\beta(v),h)$.
This uses at most $5m^2\le r$ colors.

Suppose that a monochromatic weak
$\alpha\ominus\beta$-wave
$x_1<\cdots<x_{a+b+1}$ exists.  Write $(u_i,h_i,v_i)$ for the coordinates of
$x_i$.  Monochromaticity forces a common slab coordinate, and write $h_i=h$ for
all $i$.  Let $t$ be the position of the value $1$ in $\alpha$, so the
$t$th gap is the smallest among the first $a$ gaps.

If $u_t=u_{t+1}$, then the $t$th gap equals
$v_{t+1}-v_t$ and is smaller than $L_\beta$.  In the skew sum
$\alpha\ominus\beta$, the $t$th $\alpha$-gap has rank $b+1$, so the weak-wave
inequalities force every gap in the $\beta$ part to be no larger than this
one.  On the other hand, because the slab coordinate is fixed, every change
of the coarse coordinate creates a gap at least
$5L_\beta-(L_\beta-1)>4L_\beta$.  It follows that no $\beta$-gap changes the
coarse coordinate.  The final $b+1$ points have a common coarse coordinate;
their fine coordinates are increasing, have a common $c_\beta$-color, and
inherit the weak $\beta$ inequalities.  They form a monochromatic weak
$\beta$-wave, a contradiction.

If $u_t<u_{t+1}$, then the $t$th gap is at least
$5L_\beta-(L_\beta-1)>4L_\beta$.  It is the smallest of the first $a$ gaps,
so every one of those gaps is larger than $4L_\beta$ and consequently changes
the coarse coordinate.  Thus $u_1<\cdots<u_{a+1}$.  Put
$\Delta_i=u_{i+1}-u_i$.  If
$\alpha(i)>\alpha(j)$ but $\Delta_i<\Delta_j$, then $x_{i+1}-x_i
 \le 5L_\beta\Delta_i+(L_\beta-1)
 <5L_\beta\Delta_j-(L_\beta-1)
 \le x_{j+1}-x_j$,
contrary to the weak-wave inequalities.  Therefore
$u_1<\cdots<u_{a+1}$ is a monochromatic weak $\alpha$-wave, again a
contradiction.

Hence the coloring contains no monochromatic weak $\alpha\ominus\beta$-wave. Since
$L_\alpha=\Omega(m(\log m)^p)$ and
$L_\beta=\Omega(m(\log m)^q)$, its length is $5L_\alpha L_\beta=\Omega(r(\log r)^{p+q})$,
which proves the claim.
\end{proof}

Combining the product construction with endpoint extension yields the following two oriented singleton rules that will be used in the finite classification:
$\Pw(1\ominus\beta,r)=\Omega(r(\log r)^{q+1})$ and
 $\Pw(\alpha\oplus1,r)=\Omega(r(\log r)^{p+1})$, where \(q\) and \(p\) are the corresponding exponents of \(\beta\) and \(\alpha\).

\section{ Finite classification}\label{sec:classification}
In this section, we first prove Theorems~\ref{thm:small} and~\ref{thm:S4},
which classify the patterns of length at most four. We then combine deletion
bounds with weak-wave constructions to partition $S_5$ into five classes
and reduce the length-five case to the eight exceptional alternating
patterns. The remaining upper bound for these patterns is proved in
Section~\ref{sec:completion}.
\begin{proof}[\textbf{Proof of \cref{thm:small}}]  
For $\pi=1$, every two-point subset of $[n]$ is a $1$-wave.  Hence a
$1$-wave-free set has at most one point, and the pigeonhole principle gives
$P(1,r)=r+1$.

Every two-letter permutation has no peak.  Cong's no-peak upper bound,
Lemma \ref{lem:transfer} and Corollary \ref{cor:nopeakweak} therefore give
$g(\pi,n)=\Theta(\log n)$ and
 $P(\pi,r)=\Theta(r\log r)\quad(\pi\in S_2).$

Among the six permutations in $S_3$, the four patterns
$123,213,312,321$ have no peak.  The same argument, now with $k=3$, gives
$g(\pi,n)=\Theta((\log n)^2)$ and
 $P(\pi,r)=\Theta(r(\log r)^2)$
for these four patterns.  The remaining patterns are $132$ and $231$.
In either pattern the values $1,2$ are nonadjacent, and deleting them leaves
the one-letter pattern.  Cong's two-value deletion estimate thus gives
$g(\pi,n)=O(\log n)$.  On the other hand, the consecutive two-letter
pattern $\pi[1,2]$ has no peak, so
Lemma \ref{lem:consecutive}, Lemma~\ref{lem:transfer} and Corollary \ref{cor:nopeakweak} give
$g(\pi,n)=\Omega(\log n)$ and $P(\pi,r)=\Omega(r\log r).$
The upper bound for $P(\pi,r)$ follows from
Lemma \ref{lem:transfer} (i).
\end{proof}

\begin{proof}[\textbf{Proof of \cref{thm:S4}}]
If $\pi$ has no peak, the result is Cong's no-peak theorem~\cite[Corollary~1.11]{Cong}, together with
Lemma \ref{lem:transfer} and Corollary \ref{cor:nopeakweak}.

Suppose that $\pi$ has a peak. Cong's improved upper bound~\cite[Theorem~1.14]{Cong} gives
$g(\pi,n)=O((\log n)^2)$. Consider the two consecutive three-letter
patterns
 $\sigma=\pi[1,3]$, $\tau=\pi[2,4]$.
They cannot both have a peak: a peak in $\sigma$ would give
$\pi(2)>\pi(3)$, whereas a peak in $\tau$ would give
$\pi(3)>\pi(2)$.  Hence at least one of $\sigma,\tau$ has no peak.  By
Corollary \ref{cor:nopeakweak}, its weak-wave coloring threshold has order at least
$r(\log r)^2$.  Consecutive-pattern monotonicity and
Lemma \ref{lem:transfer} yield
$g(\pi,n)=\Omega((\log n)^2)$.  
The preceding weak-wave lower bound also gives $P(\pi,r)=\Omega(r(\log r)^2)$, since $P(\pi,r)\ge\Pw(\pi,r)$. The upper bound $P(\pi,r)=O(r(\log r)^2)$ follows from $g(\pi,n)=O((\log n)^2)$ and Lemma~\ref{lem:transfer}\textup{(i)}.
\end{proof}

We now turn to patterns of length five and recall the two deletion bounds used in the reduction.
Let $\pi^-=\pi\setminus\{1\}$, where the remaining entries are normalized.
Cong~\cite{Cong} proved
$g(\pi,n)\le C_\pi(\log n)g(\pi^-,n)$.
If $1,2$ are nonadjacent, let $\pi^{--}=\pi\setminus\{1,2\}$.
Cong~\cite{Cong} proved
$g(\pi,n)\le C'_\pi(\log n)g(\pi^{--},n)$.
For $S_5$, the following lemma identifies exactly when these recurrences can
produce an exponent-$2$ upper bound.

\begin{lemma}\label{lem:E2upper}
Among the upper bounds obtained by recursively applying Cong's one-value and
nonadjacent two-value deletion inequalities, an exponent-$2$ estimate
$g(\pi,n)=O((\log n)^2)$
for $\pi\in S_5$ occurs if and only if
$1,2\text{ are nonadjacent in }\pi$ and $\pi\setminus\{1,2\}\in\{132,231\}.$
In all remaining cases, these inequalities together with Cong's no-peak and
peak bounds supply exponent $4$ when $\pi$ has no peak and exponent $3$ when
$\pi$ has a peak.
\end{lemma}

\begin{proof}
If the stated condition holds, the two-value deletion theorem first
reduces $\pi$ to $132$ or $231$.  In each of these length-three patterns the
values $1,2$ are again nonadjacent, so a second two-value deletion reduces
to the one-letter pattern. This gives exponent $2$.

Conversely, suppose that the deletion recurrences yield an exponent-$2$
upper bound. A first one-value deletion leaves a length-four pattern.  By
\cref{thm:S4}, every such pattern has extremal order at least
$\Omega((\log n)^2)$, so multiplying by the additional logarithm from the
recurrence could give no better than exponent $3$. The first step must
 be the two-value deletion, and the residual pattern
$\sigma=\pi\setminus\{1,2\}$ has length three. To obtain exponent $2$,
we have $g(\sigma,n)=O(\log n)$.  By Theorem \ref{thm:small}, this occurs exactly
for $\sigma\in\{132,231\}$. This proves the converse.  Outside this class,
Cong's no-peak theorem gives exponent $4$ when $\pi$ has no peak, and his
general peak bound gives exponent $3$ when $\pi$ has a peak.
\end{proof}

Every length-five pattern has a no-peak consecutive block of length three.
The three length-three windows are centered at positions $2,3,4$.  A window
has a peak precisely when its center is a peak of the original permutation,
and two adjacent positions cannot both be peaks. Hence the three centers
cannot all be peaks, proving the claim. Therefore
Lemma \ref{lem:consecutive} and Corollary~\ref{cor:nopeakweak} give the universal lower bound
\begin{equation}\label{eq:universal-lower}
 g(\pi,n)=\Omega((\log n)^2),
 \qquad
 P(\pi,r)=\Omega(r(\log r)^2)
 \quad(\pi\in S_5).
\end{equation}
This already matches the upper bound for the 24 patterns in $\cE_2$.

For the exponent-$3$ lower bounds, let $\cN_5$ denote the no-peak class and
define the \emph{end-block class} $\cB_3$ to consist of the patterns
$\pi\in S_5\setminus(\cN_5\cup\cE_2)$ for which at least one of
$\pi[1,4]$ and $\pi[2,5]$ has no peak. Consecutive-pattern monotonicity and Theorem \ref{thm:S4} immediately give
an exponent-$3$ lower bound for every $\pi\in\cB_3$.

 Direct and skew sums give lower bounds for another class of patterns.  Put
\begin{align*}
\cD_1={}&\{\alpha\op\beta:
 \alpha\in\{12,21\},\ \beta\in\{213,312,321\}\},\\
\cD_2={}&\{\alpha\ominus\beta:
 \alpha\in\{123,213\},\ \beta\in\{12,21\}\},\\
\cD_3={}&\{\gamma\oplus1:
 \gamma\in\{1342,2341,3142,3241,4132,4231\}\},\\
\cD_4={}&\{1\ominus\delta:
 \delta\in\{1324,1423,1432,2314,2413,2431,3412,3421\}\},
\end{align*}
and call
$\cD_{\mathrm{prod}}=\cD_1\cup\cD_2\cup\cD_3\cup\cD_4$ the
\emph{product-certificate class}. The notation distinguishes this union
from its third subfamily $\cD_3$. These four families are disjoint and have
total size $6+4+6+8=24$.

\begin{lemma}\label{lem:Dlower}
For every $\pi\in\cD_{\mathrm{prod}}$, we have 
$\Pw(\pi,r)=\Omega(r(\log r)^3).$
Consequently,
$g(\pi,n)=\Omega((\log n)^3)$.
\end{lemma}

\begin{proof}
Every two-letter pattern has weak-wave lower bound
$\Omega(r\log r)$, and every no-peak three-letter pattern has weak-wave
lower bound $\Omega(r(\log r)^2)$. Hence Proposition \ref{prop:product} gives the
required exponent-$3$ lower bound for $\cD_1$ and $\cD_2$.

Every four-letter component appearing in $\cD_3$ has a peak at position $3$,
and every component appearing in $\cD_4$ has a peak at position $2$.  Thus
Theorem \ref{thm:S4} gives exponent $2$; more specifically, the proof of that theorem
supplies the required weak lower bound through a no-peak consecutive
three-letter pattern.  The
terminal singleton in $\gamma\oplus1$ is the maximum value, and the initial
singleton in $1\ominus\delta$ is the maximum value.  Thus
Lemma \ref{lem:endpoint} adds one logarithm, proving exponent $3$ for
$\cD_3$ and $\cD_4$.
\end{proof}

The following lemma partitions the length-five patterns into five structural classes.

\begin{lemma}\label{lem:partition}
The five sets
$\cN_5, \cE_2, \cB_3, \cD_{\mathrm{prod}},$ $\cA_5$
form a disjoint partition of $S_5$, with corresponding cardinalities
$16,24,48,24,8.$
\end{lemma}

\begin{proof}
A  permutation with no peak is obtained recursively by placing the current maximum
at one of the two ends. This gives $|\cN_5|=2^4=16$.

For a pattern in $\cE_2$, the order of the values $3,4,5$ is either $132$ or
$231$, and the ordered positions of $1,2$ form a nonadjacent ordered pair in
$[5]$.  There are $5\cdot4-2\cdot4=12$ such ordered pairs, so
$|\cE_2|=24$.  In either allowed order the value $5$ lies, in position order,
between the other two values from $\{3,4,5\}$. Let their positions be
$p<q<r$, so the value at position $q$ is $5$. For the first four-entry end
block, if $q\in\{2,3\}$ then position $q$ itself is a peak. If $q=4$, then
$r=5$; among positions $1,2,3$, the two positions occupied by $1,2$ must be
nonadjacent, which forces $p=2$. The entry at position $2$ is then one of
$3,4$, while its two neighbours are $1,2$, so position $2$ is a peak. The
case $q=5$ is impossible because $r>q$. The reflected argument proves that
the last four-entry end block also has a peak: when $q=2$, nonadjacency
forces the later large value to position $4$, and otherwise the value $5$
is already an internal peak. Hence $\cE_2$ is disjoint from
$\cN_5\cup\cB_3$.

There are $5\cdot2^3=40$ patterns whose first four entries have no peak, and
the same number whose last four entries have no peak. Both end blocks have
no peak exactly when the full permutation has no peak: the first block
forbids peaks at positions $2,3$, and the second forbids peaks at positions
$3,4$. Inclusion--exclusion therefore gives $40+40-16=64$ patterns having
at least one no-peak end block. Removing $\cN_5$ gives
$|\cB_3|=48$.

There are $3!2!=12$ patterns in which the values $1,2,3$ occupy the odd
positions and the values $4,5$ occupy the even positions. Exactly four of them lie in $\cE_2$,
$14253,24153,35142,35241,$
and the other eight are precisely $\cA_5$.  Every member of $\cA_5$ has a
peak in each end block, so it is disjoint from $\cN_5\cup\cB_3$; it is
disjoint from $\cE_2$ by construction.

Expanding the four sum families defining $\cD_{\mathrm{prod}}$ gives the following
$24$ distinct permutations:
\begin{align*}
&12435,12534,12543,13425,21435,21534,21543,23415,\\
&31425,32415,34512,34521,41325,42315,43512,43521,\\
&51324,51423,51432,52314,52413,52431,53412,53421.
\end{align*}
Each word in this display satisfies
 $\pi(3)>\pi(2),\pi(4)$, and thus position $3$ is a peak in both four-entry end blocks. Thus none belongs
to $\cN_5\cup\cB_3$, and the same inequality also excludes $\cA_5$, whose middle position
is a valley. It remains to show that
$\cD_{\mathrm{prod}}\cap\cE_2=\varnothing$.
Among the $24$ permutations listed above, exactly $12$ have
the values $1$ and $2$ in adjacent positions and therefore
do not belong to $\cE_2$.
The remaining $12$ permutations are grouped in Table~\ref{tab:deletion-patterns}
according to the normalized pattern obtained by deleting these two values.
\begin{table}[!htbp]
\centering
\caption{The twelve permutations in $\cD_{\mathrm{prod}}$ with
nonadjacent $1$ and $2$.}\label{tab:deletion-patterns}
\label{tab:corresponding_permutations}
\smallskip
\begin{tabular}{@{}cl@{}}
\toprule
  \(\pi\setminus\{1,2\}\) & Corresponding permutations \(\pi\) \\
\midrule
\(123\) & \(13425,\ 23415,\ 31425,\ 32415\) \\
\(213\) & \(41325,\ 42315\) \\
\(312\) & \(51324,\ 52314\) \\
\(321\) & \(51423,\ 51432,\ 52413,\ 52431\) \\
\bottomrule
\end{tabular}
\end{table}

None of the above patterns is $132$ or $231$.
Thus none of these permutations belongs to $\cE_2$, and hence
$\cD_{\mathrm{prod}}\cap\cE_2=\varnothing$. The five sets are pairwise disjoint, and the five cardinalities add to
$16+24+48+24+8=120=|S_5|$, so the five sets form a partition.
\end{proof}
The finite classification leaves a single task. The bounds
for $\cN_5$, $\cE_2$, $\cB_3$, and $\cD_{\mathrm{prod}}$ already match,
and \eqref{eq:universal-lower} gives the required lower bounds for every
pattern in $\cA_5$. Thus, to complete \cref{thm:S5}, it remains to prove
$g(\pi,n)=O((\log n)^2)$ for $\pi\in\cA_5$; the corresponding upper
bound for $P(\pi,r)$ then follows from Lemma \ref{lem:transfer}.
\section{Separator trees and bounds for labelled binary trees}\label{sec:separator}
We develop two tree-theoretic tools for the exceptional patterns. First,
a wave-free set gives a rooted ordered binary tree whose associated
point-block diameters decrease geometrically along each branch. We then
bound the order of a labelled binary tree in terms of its number of labels
and the overlaps between the label sets of its child subtrees.
\subsection{Heavy-separator trees}\label{sec:separator1}
 For nonempty finite sets $B,C\subset\mathbb R$,
we write $B<C$ when $\max B<\min C$, and set
$\operatorname{diam}(B)=\max B-\min B$.
Let $A=\{x_0<x_1<\cdots<x_m\}$ with $m\ge1$.  Write $a_i=x_i-x_{i-1} $ for $ 1\le i\le m$; the sequence $(a_{1},\ldots,a_{m})$ is the adjacent-gap word of $A$.

\begin{proposition}\label{prop:block-sum}
Let $\pi\in S_k$.  The set $A$ contains a $\pi$-wave if and only if there
are indices
$0\le i_0<i_1<\cdots<i_k\le m$
such that the $k$ consecutive block sums
$S_j=\sum_{t=i_{j-1}+1}^{i_j}a_t$
have relative order $\pi$, for $ 1\le j\le k$.
\end{proposition}

\begin{proof}
The selected points are $x_{i_0}<\cdots<x_{i_k}$, and their $j$-th
consecutive gap is
$x_{i_j}-x_{i_{j-1}}
 =\sum_{t=i_{j-1}+1}^{i_j}(x_t-x_{t-1})=S_j.$
The equivalence follows.
\end{proof}

\begin{lemma}\label{lem:balanced-blocks}
Let $\pi\in S_k$. Let
$W_k=1+3+\cdots+3^{k-1}=\frac{3^k-1}{2}$ and $ \eta_k=\frac1{8W_k}.$
Let $a_1,\dots,a_m$ be positive real numbers of total mass $T$.  If
$\max_i a_i\le \eta_kT,$
then the word contains $k$ nonempty consecutive blocks whose sums have
relative order $\pi$.
\end{lemma}

\begin{proof}
Set $\mu=T/(2W_k)$ and prescribe the target mass
$t_i=3^{\pi(i)-1}\mu$ for the $i$th block. Starting from the left, end
each block at the first entry for which the sum reaches or exceeds its target.
If $\Delta=\max_i a_i$, every completed block has sum in
$[t_i,t_i+\Delta)$.

We first verify that the procedure completes all $k$ blocks. We have
$\sum_i t_i=\mu W_k=T/2$, and $k\le W_k$. If the procedure first failed
while trying to complete block $j$, the remaining mass would be less than
$t_j$, while each earlier completed block would have mass less than
$t_i+\Delta$. Hence $T<\sum_{i<j}(t_i+\Delta)+t_j
 \le \sum_{i=1}^k t_i+k\Delta
 \le \frac T2+\frac{kT}{8W_k}
 \le \frac{5T}{8}$, a contradiction.
Finally, if $\pi(i)<\pi(j)$, then the two target masses differ by at least
$3\mu-\mu=2\mu$.  Since
$\Delta\le T/(8W_k)=\mu/4$, we have
$t_i+\Delta<t_j$. The $i$-th block sum is less than $t_i+\Delta$, whereas
the $j$-th block sum is at least $t_j$. Thus the  block sums
have relative order $\pi$.
\end{proof}

\begin{corollary}\label{cor:heavy-interval}
Let $A\subseteq[n]$ be $\pi$-wave-free. Every nonempty consecutive subword
$a_u,\dots,a_v$ satisfies
\[
 \max_{u\le i\le v}a_i>
 \eta_k\sum_{i=u}^v a_i.
\]
\end{corollary}

\begin{proof}
Otherwise Lemma \ref{lem:balanced-blocks}, applied to that subword, and then
Proposition \ref{prop:block-sum} would produce a $\pi$-wave in $A$.
\end{proof}

Fix $\pi\in S_k$ and assume for the remainder of this subsection that $A\subseteq[n]$ is $\pi$-wave-free and $|A|\ge2$. 
The \emph{heavy-separator tree} $T_A$ is the rooted ordered binary tree obtained recursively by selecting the leftmost largest entry of each nonempty consecutive subword as its root and recursing on the subwords to its left and right. 
 If a vertex
$v$ corresponds to the subword
$W_v=(a_p,a_{p+1},\dots,a_q)$, let $h(v)$ be the index of the selected
leftmost largest entry and define the associated consecutive point block and
separator pair by
$J_v=\{x_{p-1},x_p,\dots,x_q\},
 P_v=\{x_{h(v)-1},x_{h(v)}\}$.
For any positive length $\ell$, its \emph{dyadic scale} is
$\floor{\log_2\ell}$; a pair is said to have scale $t$ when its pair gap has
that scale.  Let
$M(v)=\sum_{i=p}^q a_i$,$
  \lambda(v)=a_{h(v)}=\max_{p\le i\le q}a_i$, and $s(v)=\floor{\log_2\lambda(v)}$.
Let $d(v)$ be its depth, measured in edges from the root.  Thus
$\operatorname{diam}(J_v)=M(v)$, the pair gap of $P_v$ is the separator
gap at $v$ and has length $\lambda(v)$, and $s(v)$ is the dyadic scale of
$v$.  Every original gap corresponds to exactly one vertex, and hence
$|V(T_A)|=m=|A|-1$.

\begin{corollary}\label{cor:heavy-tree}
Let $A\subseteq[n]$ be $\pi$-wave-free, where $\pi\in S_k$, and suppose $|A|\ge2$. The height $H=\max_{v\in V(T_A)}d(v)$ of $T_A$ satisfies
\[
 H\le 1+\frac{\log n}{-\log(1-\eta_k)}=O_\pi(\log n).
\]
\end{corollary}

\begin{proof}
Let $v_0,v_1,\dots,v_h$ be a root-to-vertex path. If the point block associated with a vertex has
mass $T$, its separator has mass greater than $\eta_kT$ by
Corollary \ref{cor:heavy-interval}, so each child has mass less than
$(1-\eta_k)T$. Since the root mass is at most $n-1<n$, iteration gives
 $1\le M(v_h)<n(1-\eta_k)^h$.
Taking logarithms yields
$h<\log n/[-\log(1-\eta_k)]$. Maximizing over all vertices gives the stated bound.
\end{proof}

For the remainder, assume $k=5$ and let $\eta_*=\eta_5$ and fix
$K=\min\{j\ge1:(1-\eta_*)^j<\eta_*/8\}$.
 This constant depends only  on the fixed pattern length.

\begin{lemma}\label{lem:bounded-lag-descent}
Let $A\subseteq [n]$ be $\pi$-wave-free for some $\pi\in S_{5}$, and let $u$ be a descendant of $v$ in $T_A$. 
If $d(u)-d(v)\ge K$, then
$M(u)<\frac{\lambda(v)}8,$
 $\lambda(u)<\frac{\lambda(v)}8,$
 and $s(u)\le s(v)-3.$
\end{lemma}

\begin{proof}
Every child point block has mass less than $(1-\eta_*)$ times the mass of its
parent. Hence
$M(u)<(1-\eta_*)^KM(v)<\frac{\eta_*}{8}M(v)$.
The separator at $v$ satisfies $\lambda(v)>\eta_*M(v)$, while
$\lambda(u)\le M(u)$, so the first two inequalities follow.  If
$s=s(v)$, then $\lambda(v)<2^{s+1}$ and hence
 $\lambda(u)<\lambda(v)/8<2^{s-2}$.
Since $s(u)$ is an integer, this implies $s(u)\le s-3$.
\end{proof}

Tree distance between two vertices is the number of edges in their unique connecting path.  A
vertex in a child subtree of $v$ is \emph{$v$-deep} if its distance from $v$
is at least $K$. A dyadic scale $t$ is \emph{bi-deep at $v$} if each child
subtree contains a $v$-deep vertex whose dyadic scale is $t$. Thus the scale occurs 
at distance at least $K$ below $v$ in both child subtrees.

\subsection{Bounds for labelled binary trees}\label{sec:labelled-trees}
A logarithmic depth bound alone does not give a quadratic
logarithmic bound on the number of vertices. To obtain that bound, we
will label the separator-tree vertices by $(s(v),d(v)\bmod K)$ and control
the labels shared by the two child subtrees. The results in this subsection
apply to arbitrary labelled binary trees and convert that overlap control
into a bound on the total number of vertices.

A \emph{labelled rooted binary tree} is a finite rooted ordered
binary tree $T$ whose vertices are assigned labels from a finite set. We
also call these labels \emph{classes}; the \emph{support} of a class is the
set of vertices carrying its label. The left or right child subtree of a
vertex consists of the child and all of its descendants; an absent child
contributes the empty subtree. For a vertex $v$, let $\Sigma_L(v)$ and
$\Sigma_R(v)$ be the sets of labels occurring in its left and right child
subtrees, respectively.

\begin{theorem}\label{thm:exact-overlap-tree}
Fix an integer $C\ge0$. Let $T$ be a finite rooted binary tree whose
vertices are labelled by $q$ classes. Assume that no class is repeated on
an ancestor--descendant pair and that
$|\Sigma_L(v)\cap\Sigma_R(v)|\le C$
for every vertex $v$. Let $D_C=2^C-1$ and $K_C=2^C$.  Then
\[
 |V(T)|\le F_C(q):=
 \begin{cases}
  2^q-1,&0\le q\le C,\\
  D_C+K_C(q-C),&q>C.
 \end{cases}
\]
Both bounds are sharp.
\end{theorem}

\begin{proof}
For $q\le C$, every root--leaf path has at most $q$ vertices, so the binary
height bound gives $|V(T)|\le2^q-1$.

For $q>C$, we prove by induction on $q$. Define
$\widetilde F_C(0)=0,
 \widetilde F_C(t)=D_C+K_C\max\{t-C,0\}
 \quad(t\ge1)$.
For $1\le t\le C$, the first part of the proof gives
$|V(T)|\le2^t-1\le D_C=\widetilde F_C(t)$, which gives the
required initial cases.  Now suppose $q>C$ and that the estimate holds for every smaller number
of classes.  Let $q_1,q_2$ be the numbers of classes in the two child
subtrees of the root.  The root class occurs in
neither child, and the overlap hypothesis gives
\begin{equation}\label{eq:child-class-sum}
 q_1+q_2\le q-1+C.
\end{equation}
If one child is empty, the nonempty child uses at most $q-1$ classes, so the
induction hypothesis and monotonicity of $\widetilde F_C$ give
$|V(T)|\le1+\widetilde F_C(q-1)\le\widetilde F_C(q)$.
Suppose both children are nonempty and put
$e_i=\max\{q_i-C,0\}$. If both $q_i>C$, then
$e_1+e_2=q_1+q_2-2C\le q-C-1$ by
\eqref{eq:child-class-sum}; if at most one exceeds $C$, the same inequality
follows from $q_i\le q-1$. Thus
 $e_1+e_2\le q-C-1$.
The induction hypothesis applied to the two child subtrees therefore gives
 $|V(T)|
 \le1+\widetilde F_C(q_1)+\widetilde F_C(q_2)
 \le1+2D_C+K_C(q-C-1)
 =D_C+K_C(q-C)$,
 because $K_C=D_C+1$.
To show that the bound is sharp when $q\le C$, take the complete binary tree with $q$ levels
and assign one class to each depth. No root--leaf path repeats a class, and
at a vertex on level $j$ the two child subtrees share at most the remaining
$q-j$ classes, which is at most $C$.

For $q>C$, let $B_C$ be the complete $C$-level tree with one of the classes
$1,\dots,C$ on each depth; for $C=0$, let $B_C$ be empty.  Set $T_C=B_C$ and,
for $q>C$, form $T_q$ from a new root of class $q$, a copy of
$B_C$ on the left carrying the same classes $1,\dots,C$, and $T_{q-1}$ on
the right.  The root class is absent below the root, the ancestry condition
holds recursively, and the two child class sets intersect in exactly the
$C$ base classes. Hence the overlap bound holds at the root and at every
other vertex by induction.  Moreover,
 $|V(T_q)|=1+D_C+|V(T_{q-1})|=D_C+K_C(q-C)$.
\end{proof}

For a real number $x$, write $x_{+}=$max$\{x,0\}$. For a class-labelled tree, denote
\[
 o(v)=|\Sigma_L(v)\cap\Sigma_R(v)|,
 \qquad
 E_C(T)=\sum_{v\in V(T)}(o(v)-C)_+.
\]
We call $o(v)$ the \emph{child-class overlap} at $v$ and $E_C(T)$ the
\emph{excess-overlap budget above $C$}.

\begin{lemma}\label{lem:variable-overlap-relabel}
Suppose that the vertices of $T$ are labelled by $Q$ classes and no class
is repeated on an ancestor--descendant pair. For every integer $C\ge0$, the same tree
can be relabelled with at most $Q+E_C(T)$ classes so that the
ancestor--descendant condition remains valid and every child-class overlap
is at most $C$. Consequently,
$|V(T)|\le F_C\bigl(Q+E_C(T)\bigr).$
\end{lemma}

\begin{proof}
Let $v$ satisfy $o(v)>C$, choose a class $\alpha$ occurring in both child
subtrees of $v$, and replace every occurrence of $\alpha$ in the right child
subtree by one fresh class $\alpha'$.  We claim that this decreases the overlap at $v$ by one and leaves it unchanged at every other vertex.

At $v$, the old class remains in the left child subtree but disappears from
the right one, while the new class occurs only on the right; hence the overlap
at $v$ decreases by exactly one. At a vertex inside the renamed right
subtree, every occurrence of $\alpha$ below that vertex is changed to
$\alpha'$ simultaneously. Thus the classes shared by its two child subtrees are unchanged except for this renaming.
Vertices in the left child
subtree of $v$ and vertices incomparable with $v$ are untouched.  Finally,
let $w$ be a proper ancestor of $v$. Both the retained occurrences of
$\alpha$ in the left child subtree of $v$ and all new occurrences of
$\alpha'$ lie in the single child branch of $w$ containing $v$.  Thus
$\alpha'$ cannot become common at $w$, and the presence of $\alpha$ in each
child branch of $w$ is exactly as before. The overlap at $w$ is unchanged.

The new support is a subset of the old support, and the remaining old support
is also a subset of the original one; both are therefore antichains in the
ancestry order.  Thus the operation preserves the ancestor--descendant
condition, decreases $E_C(T)$ by exactly one, and uses one new class. 
Repeating this operation eliminates all excess overlap in at most the initial $E_C(T)$ steps. If the
final number of classes is $Q'\le Q+E_C(T)$, then
\cref{thm:exact-overlap-tree} gives $|V(T)|\le F_C(Q')$. The function $F_C$
is nondecreasing, so $F_C(Q')\le F_C(Q+E_C(T))$, proving the stated bound.
\end{proof}
\section{ The eight exceptional alternating patterns}\label{sec:completion}

{For $\pi\in\cA_5$, we construct a $\pi$-wave by combining a three-gap wave with a two-point set.} Table~\ref{tab:exceptional_permutations} records the choices for each pattern. Here $\tau_\pi\in\{132,231\}$ is a three-gap pattern,
$\kappa_\pi\in\{\mathrm{mid},\mathrm{low}\}$ specifies the rank of one
additional internal gap, and $\dir(\pi)$ specifies whether the
$\tau_\pi$-wave occurs in the earlier or later block.
\begin{table}[!htbp]
\centering
\caption{Local completion data. The wave has gaps $(a,b,c)$, the pair has gap $e$, and $X$ is the cross-block gap. The listed five gaps have reduction $\pi$.}
\label{tab:exceptional_permutations}
\smallskip
\begin{tabular}{@{}ccccc@{}}
\toprule
$\pi$ & $\tau_\pi$ & $\kappa_\pi$ & $\dir(\pi)$ & Selected gaps\\
\midrule
14352 & 132 & $\mathrm{mid}$ & $\rightarrow$ & $(a,b,c,X,e)$\\
24351 & 132 & $\mathrm{low}$ & $\rightarrow$ & $(a,b,c,X,e)$\\
34251 & 231 & $\mathrm{low}$ & $\rightarrow$ & $(a,b,c,X,e)$\\
34152 & 231 & $\mathrm{mid}$ & $\rightarrow$ & $(a,b,c,X,e)$\\
15342 & 231 & $\mathrm{low}$ & $\leftarrow$ & $(e,X,a,b,c)$\\
25341 & 231 & $\mathrm{mid}$ & $\leftarrow$ & $(e,X,a,b,c)$\\
15243 & 132 & $\mathrm{low}$ & $\leftarrow$ & $(e,X,a,b,c)$\\
25143 & 132 & $\mathrm{mid}$ & $\leftarrow$ & $(e,X,a,b,c)$\\
\bottomrule
\end{tabular}
\end{table}

For a two-point set $\{y<z\}$, its \emph{pair gap} is $z-y$. For ordered
finite point sets $B<C$, say that they are \emph{externally separated} if
$\min C-\max B>
 \max\{\operatorname{diam}(B),\operatorname{diam}(C)\}$.
Thus every cross-block difference is larger than every difference internal
to either block. If a $\tau_\pi$-wave has gaps $(a,b,c)$, then $b$ is the
largest gap.  Put
$p=\min\{a,c\}, q=\max\{a,c\}$.
A pair gap $e$ has the required type if $p<e<q$ in the mid case and $e<p$
in the low case.

\begin{lemma}[two-block completion]\label{lem:two-block-completion}
Let $B<C$ be externally separated.
\begin{enumerate}[label=\textup{(\roman*)},leftmargin=2.1em]
\item If $\dir(\pi)=\rightarrow$, a $\tau_\pi$-wave in $B$ and a pair in
$C$ whose gap has the required type produce a $\pi$-wave.
\item If $\dir(\pi)=\leftarrow$, the same conclusion holds with the pair in
$B$ and the $\tau_\pi$-wave in $C$.
\end{enumerate}
\end{lemma}

\begin{proof}
Write $(a,b,c)$ for the three gaps of the
$\tau_\pi$-wave, $e$ for the pair gap, and $X$ for the single cross-block
gap. External separation implies that $X$ is larger than $a,b,c,e$, so it
has rank $5$. Comparing $e$ with the side gaps gives the reductions recorded in Table~\ref{tab:exceptional_permutations}. 
For example, in the first row $a<e<c<b<X$, so $\red(a,b,c,X,e)=14352$. The same comparison verifies the other seven rows, yielding precisely the patterns in $\cA_5$.
\end{proof}

{We use the completion rule to identify scales that can be counted uniformly across the separator tree.}
Fix $\pi\in\cA_5$ and let $A\subseteq[n]$ be $\pi$-wave-free with $|A|\geq2$, so that $T_{A}$ is defined.  For a vertex
$v$ of $T_A$, let
$\mathcal B(v)=\{t:t\text{ is bi-deep at }v\}, b(v)=|\mathcal B(v)|.$
Write $b=b(v)$ and list $\mathcal B(v)=\{t_1<\cdots<t_b\}$.

\begin{definition}\label{def:usable-scales}
If $\kappa_\pi=\mathrm{mid}$, let
$\mathcal U_\pi(v)=\{t_i:3\le i\le b-2\}.$
If $\kappa_\pi=\mathrm{low}$, let
$\mathcal U_\pi(v)=\{t_i:1\le i\le b-4\}.$
An empty index range gives the empty set.
\end{definition}

Note that in either case of Definition \ref{def:usable-scales},
 $|\mathcal U_\pi(v)|=(b(v)-4)_+.$

\begin{lemma}\label{lem:fixed-scale-completion}
Let $t\in\mathcal U_\pi(v)$. {The point block $J_v$ contains a}
$\tau_\pi$-wave with the following property. If $\dir(\pi)=\rightarrow$,
then every scale-$t$ pair in a {later point block externally separated from $J_v$} completes
that wave to a $\pi$-wave. If $\dir(\pi)=\leftarrow$, the same holds for
every scale-$t$ pair in an {earlier point block externally separated from $J_v$}.
\end{lemma}

\begin{proof}
Every bi-deep scale supplies, in each child subtree of $v$, a separator pair
whose gap has that scale and whose vertex is at distance at least $K$ from
$v$.  By Lemma \ref{lem:bounded-lag-descent}, every such pair gap is smaller than
$\lambda(v)/8$. Writing $P_v=\{x_{h(v)-1},x_{h(v)}\}$, {the left child point block is contained} in $(-\infty,x_{h(v)-1}]$ and {the right child point block is contained} in $[x_{h(v)},\infty)$. Take one selected pair
$\{\ell_1<\ell_2\}$ in the left child and one pair $\{r_1<r_2\}$ in the
right child. {Their crossing gap $r_1-\ell_2$ contains the separator gap at $v$ and is therefore at least $\lambda(v)$. It exceeds both pair gaps and is the largest of the three consecutive gaps in $\ell_1<\ell_2<r_1<r_2$.}

Suppose first that $\kappa_\pi=\mathrm{mid}$ and $t=t_i$, where
$3\le i\le b-2$.  Put $x=t_{i-2}$ and $z=t_{i+2}$. If
$\tau_\pi=132$, choose a scale-$x$ pair in the left child and a scale-$z$
pair in the right child; if $\tau_\pi=231$, reverse these choices.  The
crossing gap is the largest of the resulting three gaps, so the four points
form a $\tau_\pi$-wave. If $p<q$ are its side gaps and $e$ is any pair gap
of scale $t$, then
$p<2^{x+1}, e\ge2^t, e<2^{t+1}, q\ge2^z$.
Since the dyadic scales are distinct integers and
$x=t_{i-2}<t_i=t< t_{i+2}=z$, in fact $x\le t-2$ and $z\ge t+2$. Therefore
 $2p<2^{x+2}\le2^t\le e$
and
 $2e<2^{t+2}\le2^z\le q$.
Equivalently, $2p<e<q/2$, so $e$ has the required mid rank.

{If $\kappa_\pi=\mathrm{low}$ and $t=t_i$ with $i\le b-4$, define}
$y=t_{i+2}$ and $z=t_{i+4}$ and choose the scale-$y$ and scale-$z$ pairs in
opposite children in the order prescribed by $\tau_\pi$.  They again form
the side gaps of a $\tau_\pi$-wave. If $p$ is the smaller side gap, then
 $2e<2^{t+2}\le2^y\le p$,
so $e$ has the required low rank. The conclusion now follows from
Lemma \ref{lem:two-block-completion}.
\end{proof}

{The next lemma bounds the number of vertices at which a fixed scale is usable by counting them level by level.}
\begin{lemma}\label{lem:usable-scale-packing}
Fix a dyadic scale $t$ and put
$\mathcal V_t=\{v\in V(T_A):t\in\mathcal U_\pi(v)\}.$
If $H$ is the height of $T_A$, then
$|\mathcal V_t|\le2^K(H+1).$

\end{lemma}

\begin{proof}
Fix a depth $d$.  If $d<K$, the entire level has at most $2^d\le2^K$
vertices. Assume $d\ge K$ and group the depth-$d$ vertices by their ancestor
at depth $d-K$.

Suppose $u,w\in\mathcal V_t$ lie in different groups, ordered so that
$J_u<J_w$, and let $z$ be their least common ancestor.  Their ancestors at
depth $d-K$ are different, so $z$ lies strictly above depth $d-K$. Since
depths are integral,
 $d(u)-d(z)=d(w)-d(z)\ge K+1$.
By Lemma \ref{lem:bounded-lag-descent}, we have
$\operatorname{diam}(J_u)=M(u)<\lambda(z)/8,
 \operatorname{diam}(J_w)=M(w)<\lambda(z)/8$.
{The point blocks associated with the two child subtrees} are disjoint and
ordered. Thus $J_u<J_w$, and every gap joining them contains the separator
gap at $z$ and is at least $\lambda(z)$. Hence $J_u$ and $J_w$ are
externally separated.

{Each point block contains both} the $\tau_\pi$-wave supplied by
Lemma \ref{lem:fixed-scale-completion} and a pair gap of scale $t$.  If
$\dir(\pi)=\rightarrow$, use that wave in $J_u$ and the pair in
$J_w$; if $\dir(\pi)=\leftarrow$, reverse their roles. {The fixed-scale completion lemma then gives a $\pi$-wave, a contradiction.} Thus all members
of $\mathcal V_t$ at depth $d$ belong to one group, which contains at most
$2^K$ vertices. Summing over the $H+1$ levels proves the claim.
\end{proof}

\begin{theorem}\label{thm:alternating-upper}
For every $\pi\in\cA_5$ and every $\pi$-wave-free set $A\subseteq[n]$,
$|A|=O_\pi((1+\log n)^2).$
\end{theorem}

\begin{proof} If $|A|\leq 1$, the conclusion is immediate. Assume $|A|\geq 2$.
Label every vertex $v$ of $T_A$ by
 $\chi(v)=(s(v),d(v)\bmod K)$.
There are at most
 $Q=K(1+\floor{\log_2n})$
classes. {No class is repeated on an ancestor--descendant pair:} if $u$ is a
proper descendant of $v$ with the same depth residue, then
$d(u)-d(v)\ge K$, and Lemma \ref{lem:bounded-lag-descent} gives
$s(u)\le s(v)-3$.

Let $o(v)$ be the overlap of the class sets in the two child subtrees.  Call a common class \emph{shallow} if at least one of the two child
subtrees has no occurrence of that class at distance at least $K$ from $v$.
For each shallow class, choose one such deficient child subtree.  The class
is common, so it does occur in that child, and every occurrence there has
distance in $\{1,\dots,K-1\}$ from $v$.  Charge the class to one nearest such
occurrence. {Each vertex has a single label, so distinct classes are charged to distinct vertices.} Across the two child subtrees there are at most
$S_*=2+4+\cdots+2^{K-1}=2^K-2$
vertices at these distances.  Hence at most $S_*$ common classes are shallow.

Every nonshallow common class has deep occurrences in both children.  Its
scale coordinate is therefore bi-deep at $v$, and for each bi-deep scale
there are at most $K$ possible depth residues.  Thus
 $o(v)\le S_*+Kb(v)$.
Set $C_*=S_*+4K$.  By Definition \ref{def:usable-scales},
\begin{equation}\label{eq:overlap-charge}
 (o(v)-C_*)_+\le K(b(v)-4)_+
 =K|\mathcal U_\pi(v)|.
\end{equation}

{The positive gaps between points of $[n]$ occupy at most}
$1+\floor{\log_2n}$ dyadic scales. By
Corollary \ref{cor:heavy-tree} and Lemma~\ref{lem:usable-scale-packing},
$\sum_v|\mathcal U_\pi(v)|
 =\sum_t|\mathcal V_t|
 \le(1+\floor{\log_2n})2^K(H+1)
 =O_\pi((1+\log n)^2)$.

Summing \eqref{eq:overlap-charge} gives
$E_{C_*}(T_A)=O_\pi((1+\log n)^2)$.
Apply Lemma \ref{lem:variable-overlap-relabel} and then
\cref{thm:exact-overlap-tree}.  Since $C_*$ is fixed, $F_{C_*}(q)$ is bounded
by a constant for $q\le C_*$ and is linear in $q$ for $q>C_*$.  Hence, for
all $q\ge0$, $F_{C_*}(q)=O_\pi(1+q)$, and therefore
\[
 |A|-1=|V(T_A)|
 \le F_{C_*}\bigl(Q+E_{C_*}(T_A)\bigr)
 =O_\pi((1+\log n)^2).
\]
\end{proof}
\begin{corollary}\label{thm:alternating-complete}
For every $\pi\in\cA_5$, we have
$ g(\pi,n)=\Theta((\log n)^2),$ and 
 $P(\pi,r)=\Theta(r(\log r)^2).$
\end{corollary}

\begin{proof}
{The upper bound for $g(\pi,n)$ follows from \cref{thm:alternating-upper}. The matching lower bounds follow from \eqref{eq:universal-lower},} and the upper bound for
$P(\pi,r)$ follows from Lemma \ref{lem:transfer}.
\end{proof}

\begin{proof}[Proof of \cref{thm:S5}]
{The no-peak case follows from Cong's theorem. For $\pi\in\cE_2$,}
Lemma \ref{lem:E2upper} gives the exponent-$2$ upper bound, while
\eqref{eq:universal-lower} gives the matching lower bound.

Let $\pi\in\cB_3\cup\cD_{\mathrm{prod}}$.  It has a peak, so Cong's theorem gives
the exponent-$3$ upper bound.  The lower bound follows from the definition
of $\cB_3$ and Lemma \ref{lem:Dlower}, respectively.  By Lemma \ref{lem:partition},
these are exactly the $72$ patterns in part~\textup{(iii)}.

Finally, the eight patterns in $\cA_5$ are covered by
Corollary \ref{thm:alternating-complete}. All coloring assertions follow from the
weak lower bounds and Lemma \ref{lem:transfer}.
\end{proof}

\section{Concluding remarks}\label{sec:conclusion}
We determine the asymptotic orders of $g(\pi,n)$ and $P(\pi,r)$ for every permutation of length at most five. 
The length-five classification partitions $S_5$ into classes with logarithmic exponents $2$, $3$, and $4$. 
The primary new bounds are the quadratic-logarithmic upper bounds for the eight exceptional alternating patterns, 
which establish sharpness by the corresponding lower bounds.

The main new ingredient is a packing argument on a rooted ordered binary
separator tree. Each vertex $v$ receives the label $(s(v),d(v)\bmod K)$,
and each label class is an antichain in the ancestor order. The completion
and packing lemmas bound the total excess overlap between child-subtree
label sets, summed over all vertices. For $\pi\in\cA_5$ and a
$\pi$-wave-free set $A\subseteq[n]$ with $|A|\ge2$, the relabelling lemma
and the extremal tree-order bound yield $|V(T_A)|=O_\pi((1+\log n)^2)$.
Since $|A|=|V(T_A)|+1$, the same upper bound holds for $|A|$.

A natural next problem is to determine the asymptotic orders for
permutations of length six and beyond, which is still a challenging direction for future research.

\section*{Declaration of generative AI}

The authors used ChatGPT 6 Astra (OpenAI) solely for language editing and polishing of this manuscript. All revisions were reviewed and approved by the authors, who remain fully responsible for the final manuscript.

\end{document}